\documentclass[a4paper]{amsart}
\usepackage{txfonts, amsmath,amstext,amsthm,amscd,amsopn,verbatim,amssymb, amsfonts}
\usepackage{fullpage}

\usepackage[bbgreekl]{mathbbol}
\usepackage{tikz}
\usetikzlibrary{matrix}
\usetikzlibrary{shapes}
\usetikzlibrary{arrows}
\usetikzlibrary{calc,3d}
\usetikzlibrary{decorations,decorations.pathmorphing}
\usetikzlibrary{through}
\tikzset{ext/.style={circle, draw,inner sep=1pt},int/.style={circle,draw,fill,inner sep=1pt},nil/.style={inner sep=1pt}}
\tikzset{exte/.style={circle, draw,inner sep=3pt},inte/.style={circle,draw,fill,inner sep=3pt}}
\tikzset{diagram/.style={matrix of math nodes, row sep=3em, column sep=2.5em, text height=1.5ex, text depth=0.25ex}}
\tikzset{diagram2/.style={matrix of math nodes, row sep=0.5em, column sep=0.5em, text height=1.5ex, text depth=0.25ex}}
\theoremstyle{plain}

\theoremstyle{definition}

\newcommand{\R}{{\mathbb{R}}}

\newcommand{\bpm}{\begin{pmatrix}}
\newcommand{\epm}{\end{pmatrix}}

\begin{document}
\title{The obstruction to the existence of a loopless star product}

\author{Thomas Willwacher}
\address{Department of Mathematics\\ University of Zurich\\ Winterthurerstrasse 190 \\ 8057 Zurich, Switzerland}
\email{thomas.willwacher@math.uzh.ch}

\thanks{The author was partially supported by the Swiss National Science Foundation, grant PDAMP2\_137151.}

\begin{abstract}
 We show that there is an obstruction to the existence of a star product defined by Kontsevich graphs without directed cycles.
\end{abstract}

\maketitle

\section{Introduction}
The theory of deformation quantization \cite{BFFLS} studies the existence and uniqueness of $\star$-products, i.~e., of associative $\R[\hbar]$ linear products on $C^\infty(M)[[\hbar]]$ (for $M$ a smooth manifold) having the form 
\[
 f\star g= fg+\hbar m_1(f,g)+\hbar^2m_2(f,g)+\cdots
\]
where the $m_j$ are bidifferential operators.
We are interested only in the case of $M=\R^d$, $d$ very large, and of $m_j$ determined from a Poisson bivector field $\pi$ by universal formulas, by which we mean Kontsevich graphs. See \cite{K1} for the definition of those graphs.
In other words, we may replace a star product by a formal series of Kontsevich graphs
\[
 a=a_0+a_1+a_2+\cdots
\]
for our purposes, where $a_0$ is fixed to be the graph 
\[
 a_0=\begin{tikzpicture}
      \draw (0,0) -- (2,0);
\node[int] at (.5,0) {};
\node[int] at (1.5,0) {};
     \end{tikzpicture} \, .
\]
The space of formal series of Kontsevich graphs forms a graded Lie algebra and the condition of associativity of the star product translates into the Maurer-Cartan equation
\begin{equation}\label{equ:MC}
 [a,a]=0.
\end{equation}

M. Kontsevich gave explicit formulas for a series of Kontsevich graphs satisfying \eqref{equ:MC} which we denote by 
\[
 a^K=a_0^K+a^K_1+a^K_2+a^K_3+a^K_4+\cdots
\]
where $a^K_j$ is a linear combination of graphs with $j$ type I vertices.
In particular $a^K_0=a_0$ and 
\tikzset{every edge/.style={draw, -triangle 60}}
\[
 a_1^K=\frac 1 2 \begin{tikzpicture}
      \draw (0,0) -- (2,0);
\node[int] (u) at (.5,0) {};
\node[int] (w) at (1.5,0) {};
\node[int] (v) at (1,1) {};
\draw (v) edge (w) edge (u);
     \end{tikzpicture}.
\]

For the detailed construction of $a^K$ we refer the reader to the original paper \cite{K1}.
Kontsevich's formal series gives rise to a star product on $\R^d$, for finite $d$, but cannot be used to construct a star product in infinite dimensions, due to the existence of directed loops in graphs occurring in $a_K$.
It was asked whether the formula can be modified in such a manner that no graphs with directed loops occur.
This question has been studied by S. Merkulov \cite{merkulovpersonalcomm}, who showed that a loopless star product cannot exist in the graded setting and by B. Shoikhet \cite{shoikhetlinfty}, who showed that a loopless formality morphism cannot exist. 

This note contains a small calculation showing that a loopless star product does not exist in the non-graded situation either.
The result is shown by G. Dito in the recent preprint \cite{Dito} as well, using an earlier explicit calculation by Penkava and Vanhaecke \cite{penkavavanhaecke}.
The benefit of our calculation is that it is a little shorter, while not using the result of \cite{penkavavanhaecke}. 

\section{Star product}
%
%

As indicated above, we call a universal star product a star product given by a formal sum of Kontsevich graphs, with two type II vertices and with all type I vertices having exactly two outgoing edges (see \cite{K1} for the notation). 
The type I vertices formally represent copies of a Poisson bivector field $\pi$.
We furthermore identify two linear combinations of graphs if they can be transformed to each other using only the (graphical version of the) Maurer-Cartan equation $[\pi,\pi]=0$ for a Poisson bivector field. 

We shall try to construct a universal star product using only graphs without oriented loops. Note that for all star products $a_0$ and $a_1$ are fixed by definition, there is no choice. $a_2$ is the same for all loopless star products and in this case uniquely determined by the MC equation. However, note that in the Kontsevich product $a_2^K$ contains a loop graph, namely
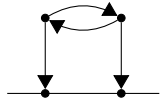
\[
\begin{tikzpicture}
\draw (-1,0)--(1,0);
\node[int] (v) at (-.5,1) {};
\node[int] (vv) at (.5,1) {};
\node[int] (w) at (-.5,0) {};
\node[int] (ww) at (.5,0) {};
\draw (v) edge[bend left] (vv) edge (w)
  (vv) edge[bend left] (v) edge (ww);
\end{tikzpicture}
\]
with nonzero weight $4\alpha$, say.
It can be removed by performing a gauge transformation using the graph
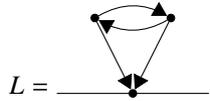
\[
L=
\begin{tikzpicture}[baseline=0]
\draw (-1,0)--(1,0);
\node[int] (v) at (-.5,1) {};
\node[int] (vv) at (.5,1) {};
\node[int] (w) at (0,0) {};
\draw (v) edge[bend left] (vv) edge (w)
  (vv) edge[bend left] (v) edge (w);
\end{tikzpicture}.
\]
So we define
\[
a := \exp(2\alpha [L,\cdot]) a^K=a_0+a_1+a_2+a_3+a_4+\cdots.
\]
Then $a_0=a_0^K$, $a_1=a_1^K$ and $a_2$ contains no graphs with directed cycles.

Now suppose our other (desirably loopless) star product reads
\[
 b = a_0+a_1+a_2+(a_3+b_3)+(a_4+b_4)+\cdots
\]
Then the Maurer-Cartan equation $[b,b]=0$ implies in particular:
\begin{align*}
 0 &= [a_0, b_3] \\
 0 &= [a_1, b_3]+[a_0, b_4]. 
\end{align*}
These imply:
\begin{align}
\label{equ:b32a}
 0 &= [a_0, b_3] \\
 0 &= [a_1, b_3]  \quad\quad\quad \text{modulo the image of $[a_0,\cdot]$}. \label{equ:b32}
\end{align}

{\bf Claim:} For any solution $b_3$ of \eqref{equ:b32a} and \eqref{equ:b32}, the sum $a_3+b_3$ contains graphs with directed cycles.

Solvability is unchanged if we add some expression $[a_0,X]$ to $b_3$, since $[a_0,a_1]=0$.
Hence we can assume that $b_3$ is in fact in the image of the (graphical version of the) Hochschild-Kostant-Rosenberg map.
In other words, it is composed of graphs for which both type II vertices have valence 1, and which are antisymmetric under interchange of these vertices.
So we may just omit the type II vertices and adjacent edges in drawings, they can be uniquely recovered.
We have the following 4 candidate graphs:

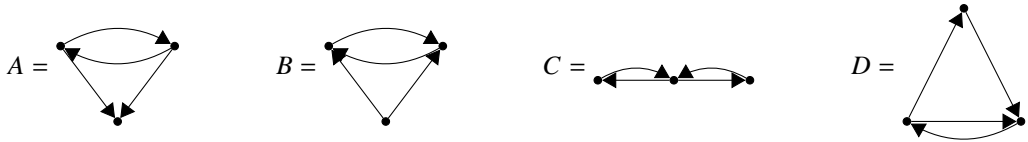
\begin{align*}
A &=
\begin{tikzpicture}[baseline={(current bounding box.center)}]
\node [int] (v1) at (-2.5,2) {};
\node [int] (v2) at (-1,2) {};
\node [int] (v3) at (-1.75,1) {};
\draw  (v1) edge[bend left] (v2);
\draw  (v2) edge (v3);
\draw  (v1) edge (v3);
\draw  (v2) edge[bend left] (v1);
\end{tikzpicture}
& 
B &=
\begin{tikzpicture}[baseline={(current bounding box.center)}]
\node [int] (v1) at (-2.5,2) {};
\node [int] (v2) at (-1,2) {};
\node [int] (v3) at (-1.75,1) {};
\draw  (v1) edge[bend left] (v2);
\draw  (v3) edge (v1);
\draw  (v3) edge (v2);
\draw  (v2) edge[bend left] (v1);
\end{tikzpicture}
&
C&=
\begin{tikzpicture}[baseline={(current bounding box.center)}]
\node [int] (v2) at (-2.5,2) {};
\node [int] (v1) at (-1.5,2) {};
\node [int] (v3) at (-0.5,2) {};
\draw  (v1) edge (v2);
\draw  (v1) edge (v3);
\draw  (v3) edge[bend right] (v1);
\draw  (v2) edge[bend left] (v1);
\end{tikzpicture}
&
D&=
\begin{tikzpicture}[baseline={(current bounding box.center)}]
\node [int] (v2) at (-1.25,3.5) {};
\node [int] (v1) at (-2,2) {};
\node [int] (v3) at (-0.5,2) {};
\draw  (v1) edge (v2);
\draw  (v1) edge (v3);
\draw  (v3) edge[bend left] (v1);
\draw  (v2) edge(v3);
\end{tikzpicture}
\end{align*}

Of these, graphs $A$, $C$ and $D$ have weight 0 in the Kontsevich star product $a^K$. Graph $B$ has non-zero weight $-\beta$, say (see \cite{vdb}, \cite{twmodules} for a computation of $\beta$). The graph $A$ attains non-zero weight $-\alpha$ in $a_3$ due to the gauge transformation, while the weights of the other graphs remain unchanged. 
 It follows that in order for $a_3+b_3$ to be loopless, $b_3$ must be a linear combination of graphs $A$ and $B$, namely $b_3=\alpha A+\beta B$.
Let us insert this $b_3$ into equation \eqref{equ:b32}. We obtain

\begin{align*}
0 &=2\alpha 
\begin{tikzpicture}[baseline={(current bounding box.center)}]
\node [int] (v1) at (-2.5,2) {};
\node [int] (v2) at (-1,2) {};
\node [int] (v3) at (-1.75,1) {};
\node [int] (v4) at (0,2) {};
\draw  (v1) edge[bend left] (v2);
\draw  (v2) edge (v3);
\draw  (v1) edge (v3);
\draw  (v2) edge[bend left] (v1);
\draw  (v2) edge (v4);
\end{tikzpicture}
+ 2\alpha 
\begin{tikzpicture}[baseline={(current bounding box.center)}]
\node [int] (v1) at (-2.5,2) {};
\node [int] (v2) at (-1,2) {};
\node [int] (v3) at (-1.75,1) {};
\node [int] (v4) at (0,2) {};
\draw  (v1) edge[bend left] (v2);
\draw  (v2) edge (v3);
\draw  (v1) edge (v3);
\draw  (v2) edge[bend left] (v1);
\draw  (v4) edge (v2);
\end{tikzpicture}
 +\alpha
 \begin{tikzpicture}[baseline={(current bounding box.center)}]
\node [int] (v1) at (-2.5,2) {};
\node [int] (v2) at (-1,2) {};
\node [int] (v3) at (-1.75,1) {};
\node [int] (v4) at (-0.5,1) {};
\draw  (v1) edge[bend left] (v2);
\draw  (v2) edge (v3);
\draw  (v1) edge (v3);
\draw  (v2) edge[bend left] (v1);
\draw  (v4) edge (v3);
\end{tikzpicture}
 +\alpha
 \begin{tikzpicture}[baseline={(current bounding box.center)}]
\node [int] (v1) at (-2.5,2) {};
\node [int] (v2) at (-1,2) {};
\node [int] (v3) at (-1.75,1) {};
\node [int] (v4) at (-0.5,1) {};
\draw  (v1) edge[bend left] (v2);
\draw  (v2) edge (v3);
\draw  (v1) edge (v3);
\draw  (v2) edge[bend left] (v1);
\draw  (v3) edge (v4);
\end{tikzpicture}
\\&\quad\quad\quad
 +2\beta 
\begin{tikzpicture}[baseline={(current bounding box.center)}]
\node [int] (v1) at (-2.5,2) {};
\node [int] (v2) at (-1,2) {};
\node [int] (v3) at (-1.75,1) {};
\node [int] (v4) at (0,2) {};
\draw  (v1) edge[bend left] (v2);
\draw  (v3) edge (v1);
\draw  (v3) edge (v2);
\draw  (v2) edge[bend left] (v1);
\draw  (v2) edge (v4);
\end{tikzpicture}
+ 2\beta 
\begin{tikzpicture}[baseline={(current bounding box.center)}]
\node [int] (v1) at (-2.5,2) {};
\node [int] (v2) at (-1,2) {};
\node [int] (v3) at (-1.75,1) {};
\node [int] (v4) at (0,2) {};
\draw  (v1) edge[bend left] (v2);
\draw  (v3) edge (v1);
\draw  (v3) edge (v2);
\draw  (v2) edge[bend left] (v1);
\draw  (v4) edge (v2);
\end{tikzpicture}
+\beta
\begin{tikzpicture}[baseline={(current bounding box.center)}]
\node [int] (v1) at (-2.5,2) {};
\node [int] (v2) at (-1,2) {};
\node [int] (v3) at (-1.75,1) {};
\node [int] (v4) at (-0.5,1) {};
\draw  (v1) edge[bend left] (v2);
\draw  (v2) edge (v3);
\draw  (v1) edge (v3);
\draw  (v2) edge[bend left] (v1);
\draw  (v4) edge (v3);
\end{tikzpicture}
\\&=\pm
2\beta
\begin{tikzpicture}[baseline={(current bounding box.center)}]
\node[int] (v1) at (-1,1.5) {};
\node[int] (v2) at (-1,-0.5) {};
\node[int] (v4) at (-2,0.5) {};
\node[int] (v3) at (0,0.5) {};
\draw  (v1) edge (v2);
\draw  (v3) edge (v2);
\draw  (v3) edge (v1);
\draw  (v4) edge (v1);
\draw  (v4) edge (v2);
\end{tikzpicture}
\pm 2(\alpha+\beta)
\begin{tikzpicture}[baseline={(current bounding box.center)}]
\node[int] (v1) at (-1,1.5) {};
\node[int] (v2) at (-1,-0.5) {};
\node[int] (v4) at (-2,0.5) {};
\node[int] (v3) at (0,0.5) {};
\draw  (v1) edge (v2);
\draw  (v1) edge (v3);
\draw  (v2) edge (v3);
\draw  (v4) edge (v1);
\draw  (v4) edge (v2);
\end{tikzpicture}\, .
\end{align*}

For the last line one must use that $[\pi, \pi]=0$. 
One checks that this combination of graphs non-zero, irrespective of the values of the non-zero numbers $\alpha$, $\beta$.
This is Shoikhet's obstruction, cf. \cite{shoikhetlinfty}. 
Note again that all vertices should be understood as having two outgoing edges. If they have less in the drawing, one must add other outgoing edges to type II vertices.

\bibliographystyle{plain}
\bibliography{../biblio} 

\end{document}